\documentclass[12pt, oneside,reqno]{amsart}   	
\usepackage{geometry}
\usepackage{float}
\allowdisplaybreaks
\usepackage{graphicx}				
\usepackage{amssymb}
\usepackage{amsmath}
\usepackage[all]{xy}
\usepackage{mathtools}
\usepackage{tikz-cd}
\usepackage{enumitem}
\usepackage{tikz-cd}
\usepackage{appendix}

\usepackage[%
bookmarks=true,
colorlinks,
linkcolor=blue,
urlcolor=blue,
citecolor=blue,
plainpages=false,
pdfpagelabels,
final,
breaklinks=true\between
]{hyperref}
\usepackage{hyperref}
\usepackage[numbers,sort&compress]{natbib}

\newtheorem{thm}{Theorem}[section]
\newtheorem{lem}[thm]{Lemma}
\newtheorem{cor}[thm]{Corollary}
\newtheorem{prop}[thm]{Proposition}

\numberwithin{equation}{section}

\def\Pb{\ifmmode{\Bbb P}\else{$\Bbb P$}\fi}
\def\Z{\ifmmode{\Bbb Z}\else{$\Bbb Z$}\fi}
\def\C{\ifmmode{\Bbb C}\else{$\Bbb C$}\fi}
\def\R{\ifmmode{\Bbb R}\else{$\Bbb R$}\fi}
\def\S{\ifmmode{S^2}\else{$S^2$}\fi}

\def\S{\cal S}

\newenvironment{pf}{\paragraph{Proof:}}{\hfill$\square$ \newline}

\begin{document}

\title[multiple ends]{An unknottedness result for self shrinkers with multiple ends} 
\author{Alexander Mramor}
\address{Department of Mathematical Science, University of Copenhagen \newline
\indent Universitetsparken 5, DK-2100 Copenhagen Ø, Denmark}
\email{almr@math.ku.dk}

\begin{abstract} In this article we prove an unknottedness result for asymptotically conical self shrinkers in $\R^3$ with multiple ends which bound a handlebody in a natural sense, using the mean curvature flow. As a corollary of this and previous work, asymptotically conical self shrinkers with two ends are unknotted.
\end{abstract}
\maketitle
\section{Introduction}

In this article we prove the following isotopy rigidity result about asymptotically conical self shrinkers with multiple ends, self shrinkers being fundamental singularity models for the mean curvature flow. Our definitions of handlebody in the noncompact case and topologically standard are given immediately below: 
\begin{thm}\label{mainthm} Let $M^2 \subset \R^3$ be an asymptotically conical self shrinker with $k \geq 2$ ends. Then if $M$ bounds a handlebody, it is topologically standard.
\end{thm} 

 As a consequence of the theorem and the author's previous work in \cite{Mra2} we show the following afterwards, of course by showing such shrinkers bound handlebodies:

\begin{cor}\label{maincor} Let $M^2 \subset \R^3$ be a properly embedded, asymptotically conical self shrinker with two ends. Then $M$ is topologically standard. 
\end{cor}

We say a closed genus $g$ surface $\Sigma \subset \R^3$ is topologically standard/unknotted when it is ambiently isotopic to a standardly embedded genus $g$ surface, by which we mean a surface given by $g$ tori embedded as tubular neighborhoods of round cricles arranged along a line with adjacent ones glued together by straight cylindrical sections when $g > 0$ and when $g = 0$ a round sphere; in particular a topologically standard torus is isotopic to the tubular neighborhood of an unknot. Abusing the definition somewhat we will often refer to the condition that $M$ bounds a region $R$ diffoemorphic to a handlebody with $k$ solid half cylinders attached simply as that $M$ bounds a handlebody, since it is a natural extension of the notion to noncompact domains, at least for our purposes. When $k > 0$ we will say such a surface $M$ is topologically standard/unknotted when it is isotopic to a standardly embedded closed surface with $k$ straight, round half cylinders attached along its convex hull -- note though for the case of multiple ends there are a number of ways one may define topological standardness depending on, say, the configuration of the ends with respect to one another as discussed more below. 
$\medskip$

Note that as stated the isotopy indicated must have unbounded speed at least to ``round out'' the ends into cylinders, although this is the only reason why it must have unbounded speed in that there are no ends being untangled from each other, or knots sent to spatial infinity, using this somehow. In particular in the course of the proof what is implied in an intermediate step is that $M$ is isotopic to a standardly embedded surface with conical ends attached along its convex hull, which are unknotted by any reasonable definition, by a finite speed isotopy (note in this dimension, the links of such ends are simply collections of circles which are all isotopic to round ones). These ends may then be rounded out easily enough as discussed more at the end of the proof. This statement above seems more satisfying given the type of result it is, though, which is why we give it instead.  
$\medskip$

That the ends are asymptotically conical is a very natural assumption with L. Wang's work and Ilmanen's cylinderical end rigidity conjecture in mind -- these are described more in section 2. Essentially, in the conjectural picture the assumption that the ends are asymptotically conical in the statements above can be dropped. Previously, analogous versions of this result were established in the compact case by the author with S. Wang in \cite{MW1} and in the noncompact case with one asymptotically conical end by the author in \cite{Mra2}; the results in these cases don't require apriori that $M$ bounds a handlebody. For this reason we have $k \geq 2$ in the statement. In a result allowing even for multiple ends, Brendle showed in \cite{Bren} that noncompact self shrinkers in $\R^3$ with no topology (roughly speaking) must coincide with either the plane or round cylinder. The result in this work is an improvement in the sense that it allows for nontrivial topology, although no geometric rigidity is claimed.  
$\medskip$

There are rigourously constructed examples of noncompact self shrinkers with one end besides the plane via desingularizations of the sphere and plane by Kapouleas, Kleene, and M{\o}ller \cite{KKM}. To the author's knowledge though there are no rigorous constructions of self shrinkers explicitly with more than one end besides the cylinder: Ketover's examples constructed in \cite{Ketover} might possibly be noncompact and hence could furnish examples of shrinkers with more than one end; see also the more recent constructions of A. Sun, Z. Wang and X. Zhou \cite{SWZ}. On the other hand there is good numerical evidence of self shrinkers with two asymptotically conical ends with interesting topology going back to Chopp \cite{Cho} (see also Ilmanen \cite{INote}). The very recent work of Buzano, Nguyen, Schulz \cite{BNS} has additional numerical examples. It's not hard to see that these 2 ended examples are handlebodies and are unknotted in the sense of the theorem above.
$\medskip$

One facet of our argument that might be of interest is that on the topological side it doesn't exclusively use Waldhausen's theorem, which says in various contexts that Heegaard splitings must be topologically standard. This is for a natural reason: Hall's examples \cite{Hall} show that minimal surfaces and hence Heegaard splitings of the ball with multiple boundary components may be knotted, which correspond roughly to surfaces with multiple ends and finite topology so that some extra assumption, or observation utilizing completeness in a special way, etc. should be necessary. Dealing with this is discussed more in the next few paragraphs, starting at classical minimal surfaces: 
$\medskip$

 There are related unknottedness results for minimal surfaces in $\R^3$, the most relevent one being the result of Meeks and Yau \cite{MY} for minimal surfaces of finite topology and any number of ends. However, there additional topological facts of complete classical minimal surfaces are used, most apparently to the author the topologically parallel ends theorem which says that a classical minimal surface $M$ with $k$ ends and finite topology is isotopic to a surface asymptotic to $k$ parallel planes -- intuitively this is true because the blowdown of a minimal surface, with some growth assumptions, yields a minimal cone, possibly with multiplicity, which in the surface case must be a plane. This allows the authors to essentially reduce down to the case that $M$ has one end, for which Waldhausen's theorem can be applied; see also the discussion at the start of section 3 below. Correspondingly their notion of topological standardness is of $k$ parallel ends connected together by straight cylindrical segments, with an unknotted compact part attached as illustrated in figure 1 of their article, which seems to differ from ours at least when there are 3 or more ends. This is because it seems that generally speaking if a surface has 3 or more parallel ends it won't be a handlebody, even if topologically standard in their defintion, in our sense. One can see this even considering just three parallel planes, adjacent planes attached by cylinders.  
 $\medskip$

 To the author's knowledge whether such a parallel ends statement also holds for shrinkers is unclear and there seems to be hints both for and against it. For instance, one might imaginably desingularize the cylinder and plane to produce examples with 3 parallel ends as mentioned in \cite{XHN} in the context of the cylindrical end rigidity conjecture, suggesting perhaps that it is true (although actually doing so may be difficult if possible). On the other hand there are certainly examples where the asymptotic behavior of self shrinkers seems to differs strongly from that of minimal surfaces, for instance that their links don't seem to satisfy a variational characterization. In the case where apriori the ends are parallel it seems plausible that their arguments can be applied, though. Resolving this apparent uncertainty one way or the other seems to be a natural question for future work.
$\medskip$ 
 
  Instead, roughly speaking in our argument we use the assumption that our shrinkers $M$ bound a handlebody to give a simple criterion for the ``compact'' part of $M$ to be topologically standardly embedded, which is satisfied by using the (renormalized) mean curvature flow along with knowledge about its long term fate either from its clearing out or that the flow becomes starshaped in the appropriate sense. The proof that the criterion is valid uses Waldhausen's theorem for Heegaard splitings of $S^3$. This is then combined with the fact that shrinker mean convex flows eventually become starshaped to finish the argument, although some care is also needed to keep track of the change in isotopy type along the flow to show theorem \ref{mainthm} in section 3. The constuction of the flow and conclusions about its long term behavior implicitly use the completeness of $M$. To show the corollary in section 4, we use the ``$\pi_1$ surjectivity criterion'' for surfaces which bound handlebodies along with a simple coloring argument (for lack of a better term) specific to $k =2$ to see that shrinkers with two ends bound handlebodies in our sense, so that the main theorem may be employed. 
  $\medskip$

$\textbf{Acknowledgments:}$  The author thanks Jacob Bernstein for discussing what was known about examples of self shrinkers with him as well as some discussions related to their spatial asymptotics; this paper was started at JHU when the author was a postdoc there and he thanks them for their hospitality. He also thanks Niels Martin M{\o}ller for his interest and for being a sounding board on some earlier iterations of the argument, and Nathalie Wahl for a valuable discussion on some topological aspects of the argument. In the course of preparing the article, he was supported by CPH-GEOTOP-DNRF151 from the Danish National Research Foundation, CF21-0680 from the Carlsberg Foundation and an AMS--Simons travel grant and is grateful for their assistance.

\section{Preliminaries}

Let $X: M \to N^{n+1}$ be an embedding of $M$ realizing it as a smooth closed hypersurface of $N$, whose image by abuse of notation we also refer to as $M$. Then the mean curvature flow $M_t$ of $M$ is given by (the image of) $X: M \times [0,T) \to N^{n+1}$ satisfying the following, where $\nu$ is the outward normal: 
\begin{equation}\label{MCF equation}
\frac{dX}{dt} = \vec{H} = -H \nu, \text{ } X(M, 0) = X(M)
\end{equation}
$\medskip$

By the comprison principle, which says that two disjoint flows will stay disjoint at least when one of them is compact, one can see singularties, that is points where the curvature blows up along the flow, occur often which makes their study important. To study these singularities, one may parabolically rescale about the developing high curvature region to obtain an ancient flow defined for times $(-\infty, T]$; when the base point is fixed this is called a $\textit{tangent flow blowup}$ which, as described by Ilmanen in his preprint \cite{I} for flows of 2 dimensional surfaces, will be modeled on \textit{smooth} self shrinkers: these are surfaces satisfying the following equivalent definitions: 
 \begin{enumerate}
\item $M^n \subset \R^{n+1}$ which satisfy $H - \frac{\langle X, \nu \rangle}{2} = 0$, where $X$ is the position vector. 
\item  Minimal surfaces in the Gaussian metric $G_{ij} = e^{\frac{-|x|^2}{2n}} \delta_{ij}$.
\item Surfaces $M$ which give rise to ancient flows $M_t$ that move by dilations by setting $M_t = \sqrt{-t} M$.
\end{enumerate}

Of course, as the degenerate neckpinch of Angenent and Velasquez \cite{AV} illustrates tangent flows do not capture quite all the information about a developing singularity but they are a natural starting point. The Gaussian metric is a poorly behaved metric in many regards; it is incomplete and by the calculations in \cite{CM1} its scalar curvature at a point $x$ is given by:
\begin{equation}
R = e^{\frac{|x|^2}{2n}}\left( n+ 1 - \frac{n-1}{4n} |x|^2 \right)
\end{equation}
We see that as $|x| \to \infty$ the scalar curvature diverges, so there is no way to complete the metric. Also since $R$ is positive for $|x|$ small and negative for large $|x|$, there is no sign on sectional or Ricci curvatures. On the other hand it is $f$-Ricci positive, in the sense of Bakry and Emery with $f = -\frac{1}{2n} |x|^2$, suggesting it should satisfy many of the same properties of true Ricci positive metrics (see \cite{WW}). Indeed, this provides some idea as to why one might expect an unknottedness result for self shrinkers, because analogous unknottedness results hold in Ricci positive metrics on $S^3$. 
$\medskip$

As is well known, the second variation for formula for area shows there are no stable minimal surfaces in Ricci positive manifolds, see for instance chapter 1 of \cite{CM}. This turns out to also be true for minimal surfaces of polynomial volume growth in $\R^n$ endowed with the Gaussian metric as discussed in \cite{CM1}. To see why this is so, the Jacobi operator for the Gaussian metric is given by: 
\begin{equation}
L = \Delta + |A|^2 - \frac{1}{2} \langle X, \nabla(\cdot) \rangle + \frac{1}{2}
\end{equation} 
The extra $\frac{1}{2}$ term is essentially the reason such self shrinkers unstable in the Gaussian metric: for example owing to the constant term its clear in the compact case from this that one could simply plug in the function ``1'' to get a variation with $Lu >0$ which doesn't change sign implying the first eigenvalue is negative. In fact, every properly embedded shrinker has polynomial volume growth by Q. Ding and Y.L. Xin: 
\begin{thm}[Theorem 1.1 of \cite{DX}]\label{proper} Any complete non-compact properly immersed self-shrinker $M^n$ in $\R^{n+m}$ has Euclidean volume growth at most. \end{thm}

It is also known that the Frenkel property for self shrinkers holds under very general conditions; see \cite{CCMS}. This is the property that any two self shrinkers must intersect -- if one of the shrinkers is compact this is an easy consequence of the comparison principle, because the origin is in the limit set as $t$ aproaches zero. These ``Ricci--like'' properties of shrinkers will play an important role below. Next we discuss the asymptotic structure of self shrinkers:  
$\medskip$

A \emph{regular cone} in $\Bbb R^3$ is a surface of the form $C_\gamma=\{r\gamma\}_{r\in (0,\infty)}$ where $\gamma$ is smooth simple closed curve in $S^2$. An end of a surface $M^2\hookrightarrow \R^3$ is \emph{asymptotically conical} with asymptotic cross section $\gamma$ if $\rho M\to C_\gamma$ in the $C^2_\mathrm{loc}$ sense of graphs as $\rho\searrow 0$ restricted to that end. 
$\medskip$

Similarly we define \emph{asymptotically cylindrical} ends to be ends which are asymptotically graphs over cylinders (with some precsribed axis and diameter) which converge to that cylinder in $C^2_{loc}$ on that end. 
$\medskip$

The reason we focus on such types of ends is the following important result of L. Wang, which says that these are the only possible types of ends which may arise in the case of finite topology: 
\begin{thm}[Theorem 1.1 of \cite{Lu}]\label{Lu-ends}

If M is an end of a noncompact self-shrinker in $\R^3$ of finite topology, then either of the following holds:
\begin{enumerate}
\item $\lim_{\tau \to \infty} \tau^{-1} M = C(M)$ in $C_{loc}^\infty(\R^3 \setminus 0)$ for $C(M)$ a regular cone in $\R^3$
\item $\lim_{\tau \to \infty} \tau^{-1} (M - \tau v(M)) = \R_{v(M)} \times S^1$ in $C_{loc}^\infty(\R^3)$ for a $v(M) \in \R^3 \setminus \{0\}$
\end{enumerate}
\end{thm}
In particular, theorem \ref{Lu-ends} applies to self shrinkers which arise as the tangent flow to compact mean curvature flows. Ilmanen in \cite{I2} conjectured that there would be no cylindrical ends unless the self shrinker was itself a cylinder: this is a natural working assumption from an analytical perspective because analysis on cones is in some sense easier than on cylinders. Later, L. Wang \cite{Lu1} gave evidence that it was true by verifying it for shrinkers very quickly asymptotic to cylinders. Note that since the convergence is in $C^\infty_{loc}$, higher multiplicity convergence is ruled out and hence the number of links of the cone must be the same as the number of ends. 
$\medskip$

In our argument, we will wish to apply the mean curvature flow on certain perturbations of a given shrinker, but it might be interesting to note it is not advantageous to do so in the Gaussian metric, which is perhaps naively the most natural setting for a flow argument. This is because the metric is poorly behaved at infinity, as one sees from the calculation of its scalar curvature, which introduces some analytic difficulites for using the flow. Instead, we consider the renormalized mean curvature flow (which we'll abreviate RMCF) defined by the following equation: 

\begin{equation}\label{renorm1}
\frac{dX}{dt} =  \vec{H} + \frac{X}{2}
\end{equation}

Where here as before $X$ is the position vector on $M$. It is related to the regular mean curvature flow by the following reparameterization; this allows one to transfer many deep theorems on the MCF to the RMCF. Supposing that $M_t$ is a mean curvature flow on $[-1,T)$, $-1 < T \leq 0$ ($T = 0$ is the case for a self shrinker). Then the renormalized flow $\hat{M_\tau}$ of $M_t$ defined on $[0, -\log(-T))$ is given by 
\begin{equation}\label{param}
\hat{X}_\tau = e^{\tau/2} X_{-e^{-\tau}},\text{ } \tau = -\log{(-t)}
\end{equation}

Note up to any finite time the reparameterization is bounded and preserves many properties of the regular MCF, like the avoidance principle. This is a natural flow for us to consider because it is up to a multiplicative term the gradient flow of the Gaussian area and fixed points with respect to it are precisely self shrinkers. More precisely, writing $H_G$ for the mean curvature of a surface with respect to the Gaussian metric:
\begin{equation}\label{Hrelation}
H_G = e^{\frac{|x|^2}{8} }(H - \frac{\langle X, \nu \rangle}{2})
\end{equation}
$\medskip$

With this in mind, the author showed in his previous article \cite{Mra2} that one can then construct a weak flow (i.e. a flow allowing for singularities), which applies to pertubations of asymptotically conical noncompact self shrinkers by the first eigenfunction to the Jacobi operator produced in \cite{BW}. In comparison to \cite{Mra2} the last two items in the statement were explicitly added for the sake of clarity, they are an easy consequence of the construction. Note that employing this is an instance where completeness is used in our proof:
\begin{thm}\label{LSF} Let $M \subset \R^3$ be a 2-sided asymptotically conical surface such that $H - \frac{\langle X, \nu \rangle}{2} \geq c(1 + |X|^2)^{-\alpha}$ for some constants $c, \alpha >0$ and choice of normal, and so that $|A(p)|^2 \to 0$ for any $p \in M \cap B(p, R)^c$. Then denoting by $K$ the region bounded by $M$ whose outward normal corresponds to the choice of normal on $M$, the level set flow $M_t$ of $M$ with respect to the renormalized mean curvature flow satisfies
\begin{enumerate}
\item inward in that $K_{t_1} \subset K_{t_2}$ for any $t_1 > t_2$, considering the corresponding motion of $K$.
\item $M_t$ is the Hausdorff limit of surgery flows $S_t^k$ with initial data $M$. 
\item $M_t$ is a forced Brakke flow (with forcing term given by position vector). 
\item Up to any fixed time, the flow is smooth with curvature bounded by any positive constant outside a suitably large ball centered at the origin.
\item $M_t$ will be $\alpha$ noncollapsed with respect to the shrinker mean curvature in any bounded spacetime region for some choice of $\alpha$ (depending on choice of domain).  
\end{enumerate}
\end{thm} 
See also \cite{HS, BH, HK, HKet, JH, JH1, Lau, Hol}. In \cite{Mra2} the author then uses this result to show asymptotically conical self shrinkers with one end must be topologically standard by using a characterization of of topological standardness due to Waldhausen discussed more in the next section. In that setting, if the surface isn't topologically standard the flow must converge to a nonempty stable self shrinker of polynomial volume growth as $t \to \infty$, by work of White \cite{W}, which contradicts that there are no such stable self shrinkers (or alternately by using the Frankel property for self shrinkers, comparing against the original shrinker). In particular, the following clearing out statement is true: 
\begin{lem}\label{Wthy} Let $M_t$ be the flow defined above in theorem \ref{LSF}. Then for any bounded region $U$, there exists $T > 0$ so that $M_T \cap U = \emptyset$. 
\end{lem}

In theorem \ref{LSF}, the order of decay assumption on shrinker mean convexity in the statement above is to ensure the flow indeed stays shrinker mean convex, using section 3 of Bernstein and L. Wang \cite{BW}. While in \cite{Mra2} it was only used to study single ended self shrinkers, theorem \ref{LSF} can be applied to self shrinkers with multiple asymptotically conical ends as well using that the convergence to the asymptotic cone is with multiplicity 1. Lemma \ref{Wthy} is also valid in this setting with multiple conical ends, and will play an important role below. 
$\medskip$

\section{Proof of theorem \ref{mainthm}}

Arguing as in \cite{Mra2} by perturbing $M$ inwards and outwards in a shrinker convex way and running the flow from theorem \ref{LSF} one can see from the clearing out lemma, lemma \ref{Wthy}, that $M$ must be a Heegaard spliting, which in this case we mean that in a very large ball it bounds regions which are both handlebodies, using a characterizations of such splittings coming from covering space theory. At this step one would wish to directly invoke a Waldhausen type theorem, the most well known one being:
\begin{thm}(\cite{Wald}) Any two Heegaard splitings of genus $g$ of $S^3$ are ambiently isotopic and, in particular, are topologically standard.
\end{thm} 
Waldhausen's theorem also applies to Heegaard splitings of the ball with one boundary component, as utilized in \cite{Mra2}, and the same statement also is true for those with two boundary components which are diffeomorphic to annuli (see corollary at the bottom of page 408 in Meeks' paper \cite{Meeks}). Unfortunately though as pointed out in the introduction, generally speaking Waldhausen's theorem doesn't apply immediately when there is more than one boundary component. This is a consequence of the knotted minimal examples of P. Hall \cite{Hall} where he produces knotted minimal surfaces, which are Heegaard splitings by proposition 2 in \cite{Meeks}, in the ball with two boundary components and genus 1; this also shows that Meeks' claim is sharp. Of course Hall's examples also suggest that to show theorem \ref{mainthm} the completeness of $M$ should play a role --one should be able to reproduce his examples for shrinkers with boundary; we point out here that our use of the flow in the argument below is at least one place where this is used. 
$\medskip$

Although perhaps a less efficient route than absolutely possible, we organize showing that $M$ is topologically standard into two related parts. First we show that it is actually given as a topologically standard (compact) handlebody with some ends attached, potentially in a complicated manner. Then we use that asymptotically the flow must become starshaped after reparameterization, and hence more or less unknotted while keeping track of what happens along the flow in between to complete the proof. 

\subsection{Unknottedness of the ``compact'' part of $\mathbf{M}$}

The main point of this subsection is to show the following, note little is claimed about the situation of the ends in comparison to the statement of the main theorem:
\begin{prop}\label{propcompact} With $M$ as in the statement of theorem \ref{mainthm} bounding a handlebody $H$ of genus $g$, $M$ is isotopic to a surface $N$ which is a standardly embedded genus $g$ surface $\Sigma$ with $k$ solid half cylinders attached. 
\end{prop}

To show unknottedness for $M$, at least the nonstandardness coming from the handles of $H$, we use the following concrete lemma. Although its use is potentially unecessary, the proof given below actaully eventually uses Waldhausen's theorem for closed surfaces coupled with the fact that $H$ is apriori a handlebody. Below, by handlebody for the noncompact case we mean in the sense given in the introduction. Also, when we write that disjoint embedded curves generate $\pi_1(H)$ of course then they must not share a common basepoint so to speak of them in this fashion is technically ill--defined; more precicely what we mean is that there are homotopies in $H$ of these curves to curves that pass through a common basepoint such that those generate $\pi_1(H)$. 
\begin{lem}\label{mainlem} Let $H$ be a genus $g$ handlebody $H$ in $\R^3$. Supposing there are $g$ disjoint embedded curves on $\partial H$ which bound embedded discs with interior in $H^c$ and generate $\pi_1(H)$, then $H$ is a closed topologically standard handlebody with $k$ smoothly embedded solid half cylinders attached. 
\end{lem}
\begin{pf}

To show the lemma we must show an isotopy of $H$ to a standardly embedded genus $g$ handlebody with ends attached in some manner, as we've defined in the introduction. Since there are $g$ and hence a finite number of discs, one may then thicken them in such a way so that $H$ with them attached is smooth and embedded. Denoting the curves in the statement by $\gamma_1, \ldots \gamma_g$ and the thickening of the discs they bound by $D_1, \ldots D_g$, denote this domain with boundary by $H^{+\cup_i D_i }$. 
$\medskip$

Since the $\gamma_i \subset \partial H$ generate $\pi_1(H)$ in the sense discussed before the lemma and gluing on a disc cannot result in new homotopy classes of loops, one can see $H^{+\cup_i D_i }$ must be simply connected so it follows from a capping argument, lemma 3.5 in \cite{Hat} and the classification of orientable surfaces that each connected component of $\partial H^{+\cup_i D_i }$ must be either be a sphere or a sphere with $\ell$ embedded half cylinders attached. By an easy calculation of the Euler characteristic using a simplical approximation of $M$ for instance, one can see gluing in a disc will either lower the genus of the boundary by one or disconnect it (and not simultaneously). Because there are exactly $g$ curves/discs in the statement and $H$ has $g$ handles then there must just be one boundary component and so $\partial H^{+\cup_i D_i }$ can be written as a sphere with $k$ embedded half cylinders attached. By sliding the discs off the ends as necessary (or, by redefining the ends to start/be attached suitably far away from the origin) we may suppose the ends and discs are disjoint. Using these discs, we proceed to isotope $\partial H^{+\cup_i D_i }$ with its ends capped to a surface which is concretely easy to see is a Heegaard splitting (essentially), so that we may apply Waldhausen's theorem. Now, as is well known the/a subset of $\partial H^{+\cup_i D_i }$ we identify with a sphere after capping suitably is isotopic to a round one, and using the isotopy extension theorem in a suitably large ball where the support of this isotopy is contained we get that $\partial H^{+\cup_i D_i }$ is isotopic to a round sphere with embedded half cylinders attached. Similarly by the isotopy extension theorem, $H^{+\cup_i D_i }$ is isotopic to a round ball $B$ with some solid half cylinders attached. The thickened discs $D_i$, which are carried along under the isotopy and remain pairwise disjoint, are diffeomorphic to $D^2 \times [0,1]$ in the obvious way, and we see that by isotoping the $D_i$ inward along the part of their boundary corresponding to $S^1 \times [0,1]$ we may isotope the $D_i$ to be pairwise disjoint thickened embedded intervals $\eta_1, \ldots \eta_g$ inside $B$ with endpoints on its boundary. Drilling along them in the obvious sense, which up to isotopy is equivalent to unattaching the $D_i$ from $H$, we see that $H$ is isotopic to a round ball $B$ with $g$ holes drilled out along the $\eta_i$ along with $k$ solid ends diffeomorphic to solid half cylinders attached. 
$\medskip$

To proceed we cap off the ends in an embedded way, so that what is left is $B$ with the isotoped discs drilled out, which by a slight abuse of notation we denote by $B \setminus \{\eta_1, \ldots \eta_g \}$. From the handlebody assumption on $H$ as defined in the introduction, $B \setminus \{\eta_1, \ldots \eta_g \}$ is a handlebody of genus $g$. One point compactifying $\R^3$ to get $S^3$, the complement of $B \setminus \{\eta_1, \ldots \eta_g \}$ after this we see is a round ball (since $\partial B$ is round) with $g$ handles attached, corresponding to the $\eta_i$; in other words the complement of $B \setminus \{\eta_1, \ldots \eta_g \}$ is a handlebody of genus $g$ as well so that $\partial(B \setminus \{\eta_1, \ldots \eta_g \})$ can be seen as a Heegaard splitting of $S^3$. Waldhausen's theorem then says it is topologically standard in $S^3 \simeq \R^3 \cup \{\infty\}$, and since the isotopy can be done to avoid a neigborhood of a point it is also topologically standard as a surface in $\R^3$. The isotopy extension theorem as above lets us carry the ends and bounded domains along under this isotopy, giving the lemma.  \end{pf}

The next idea is to use the mean curvature flow to show we may apply the lemma above. The $g = 0$ case is easy as no curves and discs need to be produced to use the lemma, so we suppose in the remainder that $g \geq 1$. By flipping the normal as needed, we may then perturb $M$ outwards from the handlebody $H$ it bounds in a shrinker mean convex way as discussed in the preliminaries, using \cite{BW}, to obtain a surface $M'$ which we relabel to be $M$. We redefine $H$ accordingly as well. We may then flow it by theorem \ref{LSF}, and the flow $M_t$ of $M$ gives a corresponding monotonic increasing (in terms of set inclusion) flow $H_t$ of $H$, as indicated in item (1) of theorem \ref{LSF}; in the theorem, $K = H^c$. The point of doing this for our purposes here is that the topology of $M_t$ eventually becomes simple, as discussed in the following:

\begin{lem}\label{long} There exists $0 << T < \infty$ for which $M_T$ is smooth for which the following holds
\begin{enumerate}
\item There is a single noncompact component $M^{nc}_T$ of $M_T$ which is diffeomorphic to a sphere with $k$ half cylinders attached.
\item The compact components of $M_t$ are diffeomorphic to spheres. 
\item For other smooth times $s \geq T$, there is only one noncompact component $M^{nc}_s$ of $M_s$ and $M^{nc}_s$ and $M^{nc}_T$ are isotopic.
 
\end{enumerate}
\end{lem}
\begin{pf}

By the clearing out lemma, lemma \ref{Wthy}, that the flow is asymptotically conical, and that the flow is set monotone outward from $H$ eventually $H_t$ must not link with $H_t^c$ after some finite time $T$, by which we mean that any closed curve in $H_t$ will not link with any closed curve in its complement. Alternately this follows from the flow being eventually asymptotically starshaped after reparameterizing (this is discussed more later on). Note that times for which $M_t$ is a smooth (possibly disconnected) hypersurface are generic by \cite{CM3} using that the flow in theorem \ref{LSF} has mean convex singularities, so that without loss of generality $M_T$ is smooth. Since all components of $M_T$ must be orientable, the choice of $T$ implies that all components of $M_t$ are either compact and diffeomorphic to spheres, giving item (2), or diffeomorphic to a sphere with some number of half cylinders attached.
$\medskip$

To show (1) it remains to show that $M_T$ must have only a single compact component, which we consider next. By item (4) of Theorem \ref{LSF} we see that the number of ends of $M_t$ won't suddenly increase any noncompact component of $M_t$ must include some subset of the (flow of the) ends of $M$. Also recall that $M$ is indicated to be connected in the statement as well -- this is actually vacuous, by the Frankel theorem for shrinkers. If eventually $M_t$ has two or more noncompact components then there must be two ends of $M$ which become disconnected from each other eventually along the flow. Supposing this occured, let $t_*$ be a time shortly before the first time this happens for which the flow is smooth. Label a pair of ends which are disconnected from each other shortly after this time by $E_1$ and $E_2$. 
$\medskip$

Let $\sigma(s)$ be an embedded line in $M_{t_*}$ so that as $s \to \infty$ it was on $E_1$ and as $s \to -\infty$ it was on $E_2$. Then from the choice of ends and $t_*$ there must shortly after be a neckpinch along $\sigma_t$ for which it corresponds to the $\R$ factor of the model cylinder up to homotopy. Denoting by $\rho$ a curve corresponding to the $S^1$ factor of this neckpinch, since the flow is monotonically into $H^c$ we see that $\rho$ bounds a disc $D$ with interior in $H^c_{t^*}$ which will later be contained in $H_t$ for $t$ slightly larger than $t_*$. Ranging over all homotopy classes for such curves $\sigma$ -- there may be multiple simultaneous nekpinches ``between'' $E_1$ and $E_2$ -- the upshot is that after flowing through these neckpinches $H_t^c$ will also have at least two unbounded components. Alternately one can use that by Alexander duality (if one wishes for, say, smooth times) that the number of connected components of $\R^3 \setminus M_t$ is one more than the number of connected components of $M_t$, that similarly $r$ disjoint noncompact components of $M_t$ give $r+1$ disjoint unbounded domains of $\R^3 \setminus M_t$, and that for all times $H_t$ is connected, which one can see using item (2) of theorem \ref{LSF}, to see this. But by item (4) of theorem \ref{LSF} for a given time $H_{t}$ outside a sufficiently large ball is a union of solid half cylinders corresponding to the ends, implying in particular that for all times there is only a single unbounded component of $H^c_{t}$. This gives a contradiction showing item (1). 
$\medskip$

With item (1) in hand denote as indicated in the statement the noncompact component of $M_T$ by $M_T^{nc}$. We see in the argument above that in fact for any time $s$ $M_s$ has only one noncompact component as well, which we'll denote by $M^{nc}_s$. From item (1) for $t > T$ $M_t$ can have only disconnecting neckpinches since there are no more handles to fill in, and by an easy calculation of the Euler characteristic must only split off spheres giving item (3). \end{pf}

With these lemmas in hand, we are now ready to prove the proposition. The idea roughly is to use that the flow becomes topologically simple in the long term, as lemma \ref{long} shows, to give discs bounded by curves on $M$ for which we may apply lemma \ref{mainlem}.

\begin{proof}[Proof of proposition \ref{propcompact}] 

It's easy to see that the approximating surgery flows from item (2) of theorem \ref{LSF} will also satisfy the conclusions of lemma \ref{long}, so for the proof of the proposition by abuse of notation when we write $M_t$ in this argument we mean such a surgery flow and $H_t$ for its associated bounded domain. 
$\medskip$

Since the flow is outward from $H$, we can see that surgering $M_t$ with a pair of surgery caps is equivalent to gluing on a thickening of an embedded disc $d$ to $H_t$, $d \subset \overline{H_t^c}$, with boundary circle laying on the neck being cut. Topologically the high curvature regions discarded under the surgery flow $M_t$ during surgery times are spheres or tori which bound balls and solid tori -- see corollary 1.25 of \cite{HK} combined with item (4) of theorem \ref{LSF}. Discarding high curvature regions diffeomorphic to a sphere is equivalent to gluing a 3-ball to $H_t$ along its boundary. Discarding a high curvature region which is a torus is equivalent to attaching a solid torus $S^1 \times D^2$ to $H_t$, but we see the difference between discarding a torus in terms of its effect on $H_t$ and simply gluing in a thickening of a disc $d \subset \overline{H_t^c}$ along a cross section of that torus is just a 3-ball. From all of this because the flow in between surgery times is an isotopy, one can succesively glue thickened discs $D_j$ onto $H$ along its boundary, corresponding to $S^1$ cross sections of high curvature regions modeled locally on cylinders, up to a given time $t$ to get a domain $H^{+ \cup_j D_j}(t)$ so that the result is isotopic to $H_t$ with some 3-balls removed. Because of this we see that since $H_T$ is simply connected, $H^{+ \cup_j D_j}(T)$ must be as well.  Of course, up to a given time only finitely many discs will be added by the finiteness of surgery times.
$\medskip$

Writing $d_j$ to be the discs whoose thickening are the discs $D_j$ above we may arrange that they are pairwise disjoint, corresponding to an isotopy of $H^{+ \cup_j D_j}(t)$ where the discs are slid off each other, and so that their boundaries lay on $M$. This isotopy can be produced by considering the surgery flow in reverse. When considering the flow in reverse surgery times correspond to necks and compact components spontaneously added, not cut/deleted as when considering the flow forward. So, considering then a $S^1$ cross section $\rho^j$, embedded, of a cylinderical region corresponding to the boundary one of the discs $d_j$ and defining the embedded curves ${\rho}^j_s$ backwards in time via the flow, it can always be perturbed away from where a given neck will soon become reattached over a small time interval shortly later. Note that since we are taking $M_s$ to be a surgery flow here, it will only need to be done finitely many times in defining a given ${\rho}^j_s$. By doing the perturbations for each of the ${\rho}^j_s$ to avoid past surgery regions appropriately, we may arrange that at a given time the curves ${\rho}^j_s$ corresponding to different $d_j$ to be disjoint since all of the $\partial d_j$ (up to isotopy the initial data, working backwards) for a given surgery time are disjoint and the flow is monotone. Using that the flow is set theoretically monotone with nonzero speed anywhere ${\rho}^j_s$ above for a given disc sweeps out an embedded annulus connecting $\rho^j$ to the curve ${\rho}_s^j$ on $M_s$, so that ${\rho}^j_s$ will bound an embedded disc with interior in $H_s^c$. Similarly all the discs produced in this fashion will be disjoint because all the timeslices for annuli corresponding to different $d_j$ (that is, the ${\rho}^j_s$) are disjoint. We see from the construction that for a given $t$ and $s < t$ $H^{+ \cup_j D_j}(t)$ and the result of gluing thickenings of these discs to $H_s$ along the ${\rho}^j_s$ are isotopic up to a finite collection of 3-balls, the discprepancy due to surgery times which are before $s$. Once $s$ is less than the first surgery time, they will in fact be isotopic. So letting $s = 0$ gives that we may suppose the $d_j$ are pairwise disjoint with boundary on $M$. 
$\medskip$

Now as mentioned above gluing in a thickened disc can either disconnect the boundary of the domain or decrease it's genus by one, so we see from lemma \ref{long} and the discussion above that by the time $T$ we can pick $g$ discs $d_1, \ldots d_g$ from the discs above (possibly nonuniquely) so that the result of gluing in their thickenings $D_1, \ldots D_g$ to $H$ attached sequentially has genus $g$ less than $H$, giving that it has boundary diffeomorphic to the sphere with $k$ ends using the classification of orientable surfaces appropriately. From this $H^{+ \cup_{i=1}^g D_i}$, that is $H$ with the $D_1, \ldots D_g$ (and just them) glued on, is simply connected. It's easy to see as a consequence the boundary loops of these discs must generate $\pi_1(H)$ in the sense discussed above lemma \ref{mainlem}, so we may apply the lemma with them to give the proposition. \end{proof}

 \subsection{Isotopy change of $\mathbf{M_t}$ under the flow and its unknottedness in the long term.}
 $\medskip$
 
In this section we assume as in the previous that $M$ has been perturbed/relabeled to be shrinker mean convex, with normal choosen so that $M_t$ flows away from $H$ using the flow from Theorem \ref{LSF}. 
$\medskip$

Even knowing that $\Sigma$ from proposition \ref{propcompact} is topologically standard actually doesn't yield as much information about the isotopy class $M$ on the whole as one might expect at first glance, because there can be synergy of sorts in terms of isotopy class between each of the components when glued together. Labeling the ends of an asymptotically conical surface $E_i$ (of course, delineating between the ``compact'' part of the surface/handlebody and its ends is arbitrary, but for this example its unimportant), consider as a toy example a hypothetical knotted cylinder, given by two ends $E_1, E_2$ attached to an embedded sphere $\Sigma$ which is the tubular neighborhood of some interval. Since the ends of $M$ are asymptotically conical we can always take the $E_i$ to be isotopic to straight round cylinders by appropriately defining where they start and $\Sigma$ ends, hiding the knottedness in $\Sigma$ in a sense, while on the other hand $\Sigma$ itself is topologically standard, concretely because one could push in along where an end would be connected as figure 1 illustrates.
$\medskip$

 \begin{figure}
\centering
\includegraphics[scale = .5]{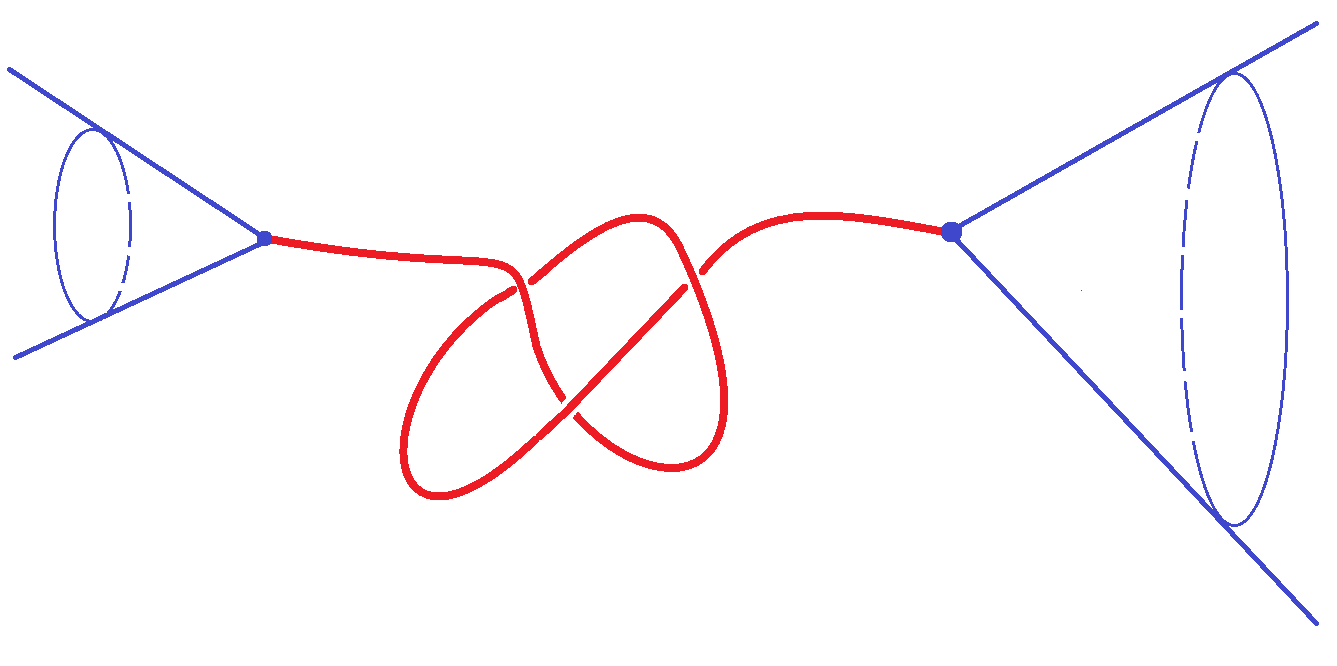}
\caption{ \small This figure shows an asymptotically conical knotted cylinder $M$ can be decomposed into three regions, two ends colored blue and a tubular neighborhood of a ``knotted interval'' diffeomorphic to a  sphere colored red, which are each individually topologically standard but combined produce a surface which isn't. Of course the blue dots represent where the ends would attach to the sphere. }
\end{figure}

To proceed, we study how the isotopy class of $M_t$ changes under the flow from $t = 0$ up to time $T$, with $T$ as in lemma \ref{long}. In the upcoming lemma and the next we use that $M$ is isotopic to the surface $N$ from proposition \ref{propcompact}. One reason for this is we know concretely how the handles of it are situated in relation to each other. This is useful, at least for psychological reasons, considering that the handles on $M$ may conceivably be arranged in a complicated way. In the statement below by $H_N$ we mean the image of $H$ under an extension of the isotopy of $M$ to $N$, and by $\overline{\gamma_i}$ we mean the longitudinal curves on the torus components of $\Sigma$ (as defined in the introduction). For reference below, one can see that by modifying the isotopy coming from Waldhausen's theorem in the proof of lemma \ref{mainlem} one can take each member of the family of surfaces to be the boundary of a round ball with some number of tubular neighborhoods of curves, some of which may have endpoints on other curves, cut out with the end result given by the boundary of a round ball drilled along $g$ straight intervals with endpoints on its boundary. Using this and that the ends are initially completely on the exterior of the ball in the proof of the lemma, by modifying the extended isotopy in the lemma along the ends appropriately (in a fixed bounded domain) to keep them on the exterior of said ball without loss of generality each of the $\gamma_i$ are nullhomotopic in $\overline{H_N^c}$. Lastly, related to the previous section by ``surger along a (embedded) disc'' below we mean performing a $1$-surgery, where the glued in discs are given by small shifts of the disc in question in a smooth, embedded way.

\begin{lem}\label{iso} The noncompact component $M^{nc}_T$ of $M_T$ is isotopic to the surface $N$ from proposition \ref{propcompact} with the handles of the genus $g$ surface $\Sigma$ filled in by surgering along a set of smoothly embedded pairwise disjoint discs $d_1, \ldots, d_g \subset \overline{H_N^c}$, where each $d_i$ has boundary ${\gamma_i} \subset N$.
\end{lem}
\begin{pf}
Recall by lemma \ref{long} there is only a single noncompact component $M_T^{nc}$ of $M_T$, $T$ as in the lemma, and the compact components of $M_T$ are diffeomorphic to spheres. We see from the proof of theorem 1.5 in the author's previous work \cite{Mra} and the approximation of the $M_t$ by surgery flows that $M$ is isotopic to $M^{nc}_T$ along with skeleta with boundary attached to $M^{nc}_T$ and some compact spherical components disjoint from $M^{nc}_T$ corresponding to the compact components of $M_T$. The isotopy in brief is constructed by ``stretching'' the necks by hand which are cut during the surgery by following the flow afterwards, while freezing the flow on high curvature regions which would normally be discarded and rounding them out appropriately. By skeleta here, we mean a thickened union of conjoined embedded intervals and embedded circles; by thickened in this lemma we mean taking the \textit{boundary} of the tubular neighborhood along the union and gluing it to $M$ in the obvious way where they meet. Denote by $S$ the union of skeleta and compact spherical components, which altogether is connected. We'll use $M_T^{nc+S}$ to refer to $M_T^{nc}$ with the skeleta and spherical components attached, isotopic to $M$ and hence $N$, and similarly define $H_T^{nc+S}$, isotopic to $H$.  Recalling that $S$ only meets $M^{nc}_T$ on the portion(s) of it corresponding to the skeleton, by an isotopy of $M_T^{nc+S}$ we may suppose that $S$ is connected and meets $M^{nc}_T$ only at a single point, by isotoping all the endpoints of the skeleton to a common point along $M^{nc}_T$. Denote the boundary circle where $S$ and $M^{nc}_T$ meet by $C$.
$\medskip$

With this in mind, as discussed preceeding the statement of the lemma without loss of generality each of the $\gamma_i$ are nullhomotopic in $\overline{H_N^c}$. So, by Dehn's lemma each bound embedded discs $d_i$ with interior in $H_N^c$ and by a cutting and pasting argument we may further arrange the $d_i$ to be pairwise disjoint. These are isotopic to discs $d_1', \ldots d_g'$ with boundary on $M_T^{nc+S}$ interior in the complement of $H_T^{nc+S}$; for reference below denote the result of surgering along these discs by $M_T^{nc+S + \cup_i d_i'}$. Also note that the result of surgering $N$ along the discs $d_i$ is isotopic to the result of surgering $M_T^{nc+S}$ along the discs $d_i'$, so that this surface is isotopic to $N$ with its handles filled in along the $d_i$. Without loss of generality $C$ lays in a plane $P$ that all the discs intersect transversely. Since the $d_i$ are all pairwise disjoint and embedded and this is preserved under isotopy we see that intersecting them with $P$ the discs correspond to disjoint simple closed curves and simple arcs with endpoints on $C$, where we note a single disc may give multiple arcs and closed curves. After isotoping the discs we may assume that they intersect $P$ only in arcs which meet $C$. 
$\medskip$

Each arc $a$ splits a disc $d_i'$ into two subdiscs, one piece $d_{i, a,S}'$ which has boundary on $S$ near $a$, and the other $d_{i, a,nc}'$ over $M_T^{nc}$ near the arc. Now, one may surger the discs along the plane $P$ by cutting them along the arcs; surgering along all arcs such discs related to $d_{i, a,nc}'$ have boundary entirely on $P \cup M_T^{nc}$, the corners of which we may smooth out, and so we may bound embedded spheres in the complement of $H_T^{nc+S}$ combining them with $P \cup M_T^{nc}$ since this union has no handles and the boundary curve on it is nullhomotopic in said complement. Note by Alexander's theorem they bound balls; in fact since there are a nonzero number of ends these spheres must be situated so that the balls they bound must lay entirely in the complement of $H_T^{nc+S}$, and we also note in the sense of set inclusion these balls come in nested families. We may thus isotope the discs across these balls to be graphs of small height over (subdomains of) $P \cup M_T^{nc}$, in particular so that they are graphs of small height over $M_T^{nc}$ away from a small collar of $C$. So, by isotopy extension and what is meant by surgery along a disc we see that $M_T^{nc+S + \cup_i d_i'}$ is isotopic to a connect sum of $M_T^{nc}$ with some surface $S'$ attached along a small shift of $C$. From the choice of the $d_i$ we get that $S'$ must be an embedded sphere attached to $M^{nc}_T$. Using Alexander's theorem again $S'$ bounds an embedded ball, so that we can isotope $S' \# M^{nc}_T$ to $M^{nc}_T$, giving the claim. \end{pf}

 From the above lemma $M^{nc}_T$ is isotopic to $N$ after filling in the handles of $\Sigma$ along the discs $d_i$. Using the organization of the handles of $N$ as an interval of adjacently attached tori, it is then easy to see that $N$ is isotopic $M^{nc}_T$ with handles added in a ``standard way,'' in that we may pick a coordinate patch of $M^{nc}_T$ and glue on $g$ small handles which bound discs and this will be isotopic to $N$ (and hence $M$). We can do this by nearly closing up the handles of $N$ along the discs $d_i$ from above; the neck of the handle of each, which can be taken to be nearly a segment of straight round cylinder since the discs are smooth, can then be slid to the edge of their handle which clearly gives a small, disc bounding handle standing ``atop'' the corresponding original handle, with now the discs $d_i$ completely filled/surgered along. Then all of them can be gathered together to be attached to $N$ within a coordinate chart as mentioned. The point of doing this is now we can write $N$ as a surface we know in a pretty concrete fashion is isotopic to $M^{nc}_T$ plus with the handles above attached, which gives then by the isotopy extension theorem that $N$ is isotopic to $M^{nc}_T$ with handles attached in a standard way on it because we can carry the handles added to $N$-with-handles-filled-in along the isotopy. This may not really be a necessary step in the proof but it seems to make using the lemma above easier to visualize -- of course in this matter there is really not quite such a thing as ``small'' since for our ultimate goal we may isotope freely. Summarizing this discussion in a lemma:

 \begin{lem}\label{added} There exists $T >> 0$ so that $N$, and hence $M$ is isotopic to $M^{nc}_t$ with $g$ handles added in a standard way for $t \geq T$. 
 \end{lem}

Our next goal is to understand $M^{nc}_T$; this is where we most strongly use that the flow is asymptotically starshaped after reparameterization: 
\begin{lem}\label{round} $M^{nc}_T$ is isotopic to a standardly embedded surface with $k$ ends and no handles.
\end{lem}
\begin{pf}

Undoing the reparameterization in the RMCF back to the ``regular'' mean curvature flow and writing it as $\widetilde{M}_t$, as observed by Bernstein and L. Wang \cite{BW} at time $t = 0$ if $\widetilde{M}_t$ is a smooth flow its must be starshaped with respect to the origin in that $x \cdot \nu \geq 0$. It was later shown that for $t$ slightly less than zero that $\widetilde{M}_t$ must be smooth even allowing for singularities by Chodosh, Choi, Mantoulidis, and Schulze \cite{CCMS}: in particular see their corollaries 8.17, 8.18 -- their flow agrees with ours since the level set flow is nonfattening. Back parameterizing note this implies for $T = T(t)$ large that $M_T$ will be smooth, after taking $T$ potentially larger. For $t$ very close to zero, denote by $\widetilde{M}^{nc}_t$ the noncompact connected component of $\widetilde{M}_t$ which corresponds to $M^{nc}_t$. Since the flow is smooth for times shortly before $0$, this gives that $\widetilde{M}^{nc}_t$ is for sufficiently small $t < 0$ isotopic to a starshaped surface (by a finite speed isotopy, using the surgery approximation), which we claim gives that $M^{nc}_T$ is a topologically standard surface with no handles by the choice of $T$ and $k$ ends by lemma \ref{long}. Reiterating the comments in the introduction such isotopy must have infinite speed to round out/cylindricalize the ends, but from how we describe it it will be clear it will not change the isotopy class of $M$ in that it would not send a knotted region to spatial infinity or ``untangle'' ends by moving the links about each other. 
$\medskip$

We describe the isotopy first for $\widetilde{M}^{nc}_t$, where $t$ is very close to zero. First isotope it to $P = \widetilde{M}^{nc}_0$, which is starshaped. Being starshaped means this is a graph over a region $D$ of a sphere centered at the origin, say $S(0,1)$ for concreteness, with some boundary components $\sigma_i$, one for each link component (in fact, they are the same curves). By the handlebody assumption for the initial data, they bound a disjoint set of discs $U_i$ which gives isotopies, via say the Riemann mapping theorem, $(\sigma_i)_s \subset U_i$ of the $\sigma_i$ to round circles. From these isotopies we get a corresponding isotopy of $P$. To see this, note that since $P$ is asymptotically a graph over the cone $C$ with link $\bigcup_i (\sigma_i)_0$ and it is starshaped, $P$ can also be defined as a graph over $D \cup (C \cap B(0, 1)^c)$, which by mollifying the corners appropriately we may assume is smooth and starshaped, so that it is isotopic to it. In the following denote by $C'_s = C_s \cap B(0, 1)^c$ where $C_s$ is the cone over $(\sigma_i)_s$, and denote by $D_s$ the corresponding region with boundary given by $(\sigma_i)_s$ so that $D_0 = D$ and $C_0 = C$. Rounding the corners of $D_s \cup C'_s$ in a continuous way in $s$ -- embeddedness is not an issue since each of the isotopies of $\sigma_i$ are contained in $U_i$, we see these furnish an isotopy from $P$ to $D_1 \cup C'_1$. Reparameterizing appropriately, at $s = 1$ each of the $(\sigma_i)_s$ are round circles so that $D_1 \cup C'_1$  is clearly graphical over a standardly embedded surface with $g = 0$ and $k$ ends in our terminology, giving the isotopy for $\widetilde{M}^{nc}_t$. By the parameterization $\widetilde{M}^{nc}_t$ is isotopic to $M^{nc}_T$, so that as claimed $M^{nc}_T$ is topologically standard. 
\end{pf}
 
\begin{proof}[Proof of Theorem \ref{mainthm}] Summarizing the above sequence of claims essentially, employing lemmas \ref{mainlem}, \ref{long} we get from proposition \ref{propcompact} that $M$ is isotopic to a surface $N$ given by a standardly embedded genus $g$ surface with $k \geq 2$ half cylinders attached, although which won't necessarily be straight and round and so could be arranged in a complicated way. By isotoping the handles of these appropriately to be ``standard handles'' added onto an aproximately planar part of $N$ in a single coordinate chart, the result which we'll call here $N'$ is isotopic to the noncompact component $M_t^{nc}$ of $M_T$, $T >> 0$ as in lemma \ref{long}, with what we called standard handles attached by lemmas \ref{iso}, \ref{added}. $M_T^{nc}$ is then isotopic to a standardly embedded surface with no handles and $k$ ends, in our sense, by lemma \ref{round} using importantly that the flow becomes starshaped as $t \to 0$ after change of parameters, so that $N'$ is isotopic to a standardly embedded surface of no handles with $g$ standard handles attached using the isotopy extension theorem. This configuration is then obviously isotopic to a standardly embedded surface with $g$ handles and $k$ ends after perhaps changing the axis direction of the cylinders somewhat, completing the argument. 
\end{proof}

\section{Proof of corollary \ref{maincor}}
In this section we show that self shrinkers with 2 asymptotically conical ends bound a handlebody in our sense, so that the work above applies.  Below by ``ends'' we mean the connected components of $M \cap B(0, r)^c$, where $r$ is taken large enough so that each component is topologically an embedded half cylinder; such an $r$ exists by the asymptotically conical assumption. Below we may similarly abuse definitions by refering to their $S^1$ cross sections as their links when the intended meaning is clear. First, we show the following which can be interpreted as saying in a sense that asymptotically such shrinkers with 2 ends bound handlebodies: 
 \begin{lem}\label{interior} Suppose $M$ is a properly embedded asymptotically conical surface whose link has two components. Then it bounds a domain $R$ with which is diffeomorphic to the interior of a domain bounded by a closed compact surface $\Sigma$ with two solid half cylinders attached. 
\end{lem} 
\begin{pf} 
 Writing the ends of $M$ as $E_1$ and $E_2$ consider the associated components of the link, $L_1$ and $L_2$; for each end there is precisely one component of the link by L. Wang's work in \cite{Lu}. Note by the Frankel property for shrinkers $M$ is connected, so since $M$ is properly embedded by Alexander duality it separates $\R^3$ into two components, a well known fact already alluded to many times before. Using this, we may color the sphere with two colors which alternate across the link of $M$, depending on which bounded component of $\R^3 \cap S(0,r)$ a point lays in. Since the link when there are 2 ends splits the sphere into 3 components the links must bound disjoint discs of the same color, giving that $M$ bounds a domain $R$ for which $B(0, r)^c \cap R$ are two disjoint solid half cylinders with boundary given by $E_1$ and $E_2$, where $r$ as in the discussion above is suitably large. Taking $\Sigma$ to be the boundary of $R \cap B(0, r)$ after smoothing the edges along $S(0,r)$, we get the statement. 
\end{pf}

Referring to the region $D \subset R$, $R$ as above, which is bounded by $\Sigma$ (capping of the ends suitably), we next wish to show $D$ is a handlebody. As discussed in \cite{L, papa}, $D$ will be a handlebody if the ``$\pi_1$ surjectivity criterion'' holds: that $\iota_*: \pi_1(\Sigma) \to \pi_1(D)$ is surjective, where $\iota$ is the inclusion map of $\Sigma$ into $D$. By covering space theory, this is equivalent to the boundary of the universal cover of $D$ being connected. Indeed, this is a criterion often applied to show a minimal surface is a Heegaard splitting, which is a surface in a 3--manifold which decomposes it into two handlebodies. 
$\medskip$

By the author's previous work \cite{Mra2}, the map $\iota_*: \pi_1(M) \to \pi_1(R)$, where $\iota^M$ is the inclusion map of $M$ into $R$, will be surjective: as a sketch this is accomplished by perturbing $M$ into $R$ in a shrinker mean convex way and flowing $M$ using Theorem \ref{LSF}. If the boundary of the universal cover of $R$ isn't connected then by an intersection number argument, considering the lift of the flow into the universal cover, one can see the flow cannot clear out giving a contradiction to lemma \ref{Wthy}. Again, by covering space theory this implies that the map above is surjective.
$\medskip$

 Note by the choice of $R$ the links of the ends $E_1$ and $E_2$ are nullhomotopic in $R$ so that $\iota_*$, a group homomorphism, restricted to $\pi_1(M)/\langle L_1, L_2 \rangle^N$ is a surjection onto $\pi_1(R)$, where we take $\langle L_1, L_2 \rangle^N$ to denote the normal closure of the subgroup of $\pi_1(M)$ generated by the links of the ends. Also, since the ends bound solid half cylinders in $R$ all closed loops in $R$ can be retracted into $D$, so that $\pi_1(R) \simeq \pi_1(D)$. Using the standard presentation of the fundamental group of a punctured surface ($M$ can be thought of as $\Sigma$ with 2 punctures), we can see that $\pi_1(M)/\langle L_1, L_2 \rangle^N \simeq \pi_1(\Sigma)$ because $\Sigma$ is obtained from $M$ by capping off the ends. Using these identifications, one can show that $\iota^\Sigma_*: \pi_1(\Sigma) \to \pi_1(D)$ will be surjective since $\iota^M_*: \pi_1(M) \to \pi_1(R)$ is. This gives that $D$ is a handlebody, so that combined with the lemma above $R$ is a handlebody in our sense as described in the introduction. With this in hand the corollary then follows from theorem \ref{mainthm}. Note this argument also gives that if an embedded self shrinker bounds a region homeomorphic to a domain bounded by a closed surface with $k \geq 2$ solid half cylinders attached, then it is topologically standard. 
$\medskip$

\end{document}